\documentclass[12pt]{article}

\usepackage{hyperref}
\usepackage{amsmath}
\usepackage{amssymb}
\usepackage{amsfonts}
\usepackage{amsopn}
\usepackage{amsthm}
\usepackage[all]{xy}

\swapnumbers
\theoremstyle{plain}
\newtheorem*{thm}{Theorem}

\newcommand{\pl}{\mathbb{P}}   
\newcommand{\zz}{\mathbb{Z}}

\DeclareMathOperator{\Spec}{\mathrm{Spec}}

\title{Essential dimension of Hermitian spaces}

\author{N.~Semenov, K.~Zainoulline\footnote{
Supported by DFG GI706/1-1 and INTAS 05-1000008-8118.\newline
}}

\date{}

\begin{document}

\maketitle

\begin{abstract}
Given an hermitian space we compute its essential dimension,
Chow motive and prove its incompressibility in certain dimensions.
\end{abstract}

The notion of an {\em essential dimension} $\dim_{es}$
is an important birational invariant of an algebraic variety 
$X$ which was introduced and studied 
by N.~Karpenko, A.~Merkurjev, Z.~Reichstein, J.-P.~Serre and others.
Roughly speaking, it is defined to be the minimal possible dimension of
a rational retract of $X$. In particular,
if it coincides with the usual dimension, then $X$ is called
{\em incompressible}.

In general, this invariant is very hard to compute. As a consequence,
there are only very few known examples of incompressible varieties:
certain Severi-Brauer varieties and projective quadrics. 
In the present paper we provide new examples
of incompressible varieties: {\em Hermitian quadrics} of dimensions
$2^r-1$.
We also give an explicit formula for the
essential dimension of a Hermitian form in the sense of O.~Izhboldin,
hence, providing a Hermitian version of the
result of Karpenko-Merkurjev \cite{KM03}. 
At the end we discuss the relations with Higher forms of Rost motives
of Vishik \cite{Vi00}.

We follow the notation of \cite{Kr07}.
Let $F$ be a base field of characteristic not $2$ and
let $L/F$ be a quadratic field extension. Let
$(W,h)$ be a non-degenerate $L/F$-Hermitian space of rank $n$
and let $q$ be the 
quadratic form associated to the hermitian form $h$ via 
$q(v)=h(v,v)$, $v\in W$.

The main objects of our study are the following
two smooth projective varieties over $F$:
\begin{itemize}
\item
the variety $V(q)$ of 
$q$-isotropic $F$-lines in $W$, i.e. a projective {\em quadric};
\item
the variety $V(h)$ of $h$-isotropic $L$-lines
in $W$ called a {\em Hermitian quadric}.
\end{itemize}

Observe that $V(q)$ has dimension $(2n-2)$ and
$V(h)$ is a $(2n-3)$-dimensional projective homogeneous variety 
under the action of the unitary group associated with $h$.
It is also a twisted form of the {\it incidence variety} that is
the variety of flags 
consisting of a dimension one and codimension one linear subspaces
in an $n$-dimensional vector space.

The forms $q$ and $h$ are closely related
by the following celebrated result
of Milnor-Husemoller (see \cite{Le79}):

\begin{quote}
A quadratic form $q$
on an $F$-vector space $W$ is the underlying form of a Hermitian form
over a quadratic field extension $L=F(\sqrt a)$ iff 
$\dim W=2n$, $q_L$ is hyperbolic,
and $\det q=(-a)^n\mod F^2$.
\end{quote}

\paragraph{1. Incompressibility}
A smooth projective $F$-variety $X$
is called {\em incompressible}
if any rational map $X\dasharrow X$ is dominant.
The basic examples of such varieties are anisotropic
quadrics of dimensions $2^r-1$ and Severi-Brauer varieties of
division algebras of prime degrees.

\begin{thm}[A]
Assume that the variety $V(h)$ is anisotropic and
$\dim V(h)=2^r-1$ for some $r>0$. Then $V(h)$ is incompressible.
\end{thm}

\begin{proof}
The key idea is that
a Hermitian quadric which is purely a geometric object can be viewed
as a twisted form of a {\em Milnor hypersurface} $M$ -- 
a topological object, namely, a generator of the Lazard ring 
of {\em algebraic cobordism} of M.~Levine and F.~Morel \cite{LM}.

More precisely, by \cite[2.5.3]{LM} the variety $M$
is the zero divisor of the line bundle  
$\mathcal{O}(1)\otimes \mathcal{O}(1)$ on $\pl^{n-1}_F\times\pl^{n-1}_F$, i.e.
it is given by the equation 
\begin{equation}\label{inceq}
\sum_{i=0}^{n-1}x_iy_i=0,
\end{equation}
where $[x_0:\ldots x_{n-1}]$ and
$[y_0:\ldots :y_{n-1}]$ are the projective coordinates
of the first and the second factor respectively.

From another hand side, the Hermitian quadric $V(h)$ is a twisted form
of the incidence variety $X=\{W_1\subset W_{n-1}\}$, where $\dim W_i=i$.
Taking $[x_0:\ldots: x_{n-1}]=W_1$ and 
$[y_0:\ldots:y_{n-1}]$ to be the normal vector
to $W_{n-1}$ we obtain that $X$ is given by the same equation \eqref{inceq},
therefore, $X\simeq M$.

By \cite[Prop.7.2]{Me02} we obtain
the following explicit formula for the Rost characteristic 
number $\eta_2$ of $M$
$$
\eta_2(M):=\frac{c_{\dim M}(-T_M)}{2}=
\frac{1}{2}{2(n-1)\choose n-1} \mod 2.
$$
It has the following property: 
\begin{equation}\label{congrel}
\eta_2(M)\equiv 1\mod 2\quad\Longleftrightarrow\quad\dim M=2^r-1\text{ for some }r>0.
\end{equation}
Since $\eta_2$ doesn't depend on the base change, 
$\eta_2(M)=\eta_2(V(h))$.

We apply now the standard arguments involving the Rost degree formula
(see \cite[\S7]{Me03}). Let $f\colon V(h)\dasharrow V(h)$ be a rational map. 
By the degree formula:
\begin{equation}\label{rostdeg}
\eta_2(V(h))\equiv \deg f\cdot\eta_2(V(h)) \mod n_{V(h)},
\end{equation}
where $n_{V(h)}$ is the greatest common divisor of degrees of all closed
points of $V(h)$. Since $V(h)$ becomes isotropic over $L$, $n_{V(h)}=2$.

Assume now that $\dim(V(h))=2^r-1$ for some $r>0$.
Then, by \eqref{congrel} $\eta_2(V(h))\equiv 1$ and by \eqref{rostdeg}
$\deg f\neq 0$ which means that $f$ is dominant. 
This finishes the proof of the theorem.
\end{proof}

\paragraph{2. Essential dimension}
Following O.~Izhboldin
we define the {\it essential dimension} of a Hermitian space $(W,h)$ as 
$$
\dim_{es}(h):=\dim V(h)-i(q)+2,
$$
where $i(q)$ stands for the first Witt index of the form $q$ 
(cf. \cite{KM03}).
The following theorem provides a {\it Hermitian version} 
of the main result of \cite{KM03}

\begin{thm}[B]
Let $Y$ be a complete $F$-variety with all closed
points of even degree. Suppose that $Y$ has a closed point of odd degree
over $F(V(h))$. Then $\dim_{es}(h)\le\dim Y$. 
Moreover, if $\dim_{es}(h)=\dim Y$,
then $V(h)$ is isotropic over $F(Y)$.
\end{thm}
\begin{proof}

In \cite{Kr07} D.~Krashen constructed
a $\pl^1$-bundle 
\begin{equation}\label{blu}
Bl_{S}(V(q)) \to V(h),
\end{equation}
where $Bl_S(V(q))$ is the blow-up of the quadric $V(q)$ 
along the linear subspace $S=\pl^{n-1}_L$.
In particular, the function field of $V(q)$ is
a purely transcendental extension of the function field of 
$V(h)$ of degree $1$,
and, therefore, our theorem follows from \cite[Theorem~3.1]{KM03}.
\end{proof}

Using Theorem (B) we can give an alternative proof of Theorem (A):
\begin{proof}[Another proof of (A)]
Let $Y$ be the closure of the image of a rational map $V(h)\dasharrow V(h)$.
Then by Theorem (B) the incompressibility of $V(h)$ follows from the equality
$\dim_{es}(h)=\dim V(h)$. The latter can be deduced from the
Hoffmann's conjecture (proven in \cite{Ka03}) if $V(h)$ is anisotropic
and $\dim V(h)=2^r-1$. Indeed, if $\dim V(h)=2^r-1$, then $\dim q=2^r+2$.
Therefore, $i(q)=1$ or $2$. But by the result of Milnor-Husemoller 
$i(q)$ must be even. Hence, $\dim_{es}(h)=\dim V(h)$.
\end{proof}

\paragraph{3. Chow motives}
We follow the notation of \cite[\S6]{CM06}.
As a direct consequence of the fibration \eqref{blu} and
the Krull-Schmidt theorem proven in \cite{CM06} we obtain the following
expressions for the Chow motives of $V(q)$ and $V(h)$:

\begin{thm}[C]
There exists a motive $N_h$ such that
\begin{equation}\label{mquad}
M(V(q))\simeq 
\begin{cases}
N_h\oplus N_h\{1\}, &\text{ if }n\text{ is even};\\
N_h\oplus M(\Spec L)\{n-1\}\oplus N_h\{1\}, &\text{ if }n\text{ is odd};
\end{cases}
\end{equation}
and
\begin{equation}\label{mvh}
M(V(h))\simeq 
\begin{cases}
N_h\oplus \bigoplus_{i=0}^{(n-4)/2}M(\pl_L^{n-1})\{2i+1\}, &\text{ if }n\text{ is even};\\
N_h\oplus \bigoplus_{i=0}^{(n-3)/2}M(\pl_L^{n-2})\{2i+1\},&\text{ if }n\text{ is odd}.
\end{cases}
\end{equation}
Observe that by the projective bundle theorem
$M(\pl_L^m)\simeq\bigoplus_{i=0}^mM(\mathrm{Spec}\, L)\{i\}$.
\end{thm}

\begin{proof}
Using the $\pl^1$-fibration \eqref{blu} 
D.~Krashen provided the following formula
relating the Chow motives of $V(q)$ and $V(h)$:
\begin{equation}\label{kreq}
M(V(q))\oplus \bigoplus_{i=1}^{n-2} M(\pl_L^{n-1})\{i\} \simeq 
M(V(h))\oplus M(V(h))\{1\}.
\end{equation}

Observe that
the motives of all varieties participating in the
formula \eqref{kreq} 
split over $L$
into direct sums of twisted Tate motives $\zz_L$.
For each such decomposition 
$M_L\simeq \bigoplus_{i\ge 0} \zz_L\{i\}^{\oplus a_i}$
we define the respective Poincar\'e polynomial by
$P_{M_L}(t):=\sum_{i\ge 0} a_it^i$.
Using the standard combinatorial description of the cellular structure on
$V(q)_L$, $V(h)_L$ and $\pl_L^{n-1}$ (see \cite{CM06}) we obtain
the following explicit formulae:
\begin{equation}\label{genf}
P_{V(q)_L}(t)=\tfrac{(1-t^n)(1+t^{n-1})}{1-t},\;
P_{V(h)_L}(t)=\tfrac{(1-t^n)(1-t^{n-1})}{(1-t)^2}\text{ and }
P_{\Spec L}(t)=2.
\end{equation}

Consider the subcategory of the category of Chow motives with
$\zz/2$-coefficients generated by $M(V(h);\zz/2)$.
Since $V(h)$ is a projective homogeneous variety,
the Krull-Schmidt theorem and
the cancellation theorem hold in this subcategory by \cite[Cor.35]{CM06}.
In particular, two decompositions of the formula 
\eqref{kreq} have to consist from the same indecomposable
summands. 

Analyzing their Poincar\'e polynomials over $L$ 
using \eqref{genf} we obtain the formulae \eqref{mquad} and
\eqref{mvh} for motives with $\zz/2$-coefficients.
Finally, applying \cite[Thm.2.16]{PSZ} 
for $m=2$ we obtain the desired formulae integrally.
\end{proof}

\paragraph{4. Higher forms of Rost motives}
In \cite[Thm.5.1]{Vi00} 
A.~Vishik proved that 
given a quadratic form $q$ over $F$
divisible by an $m$-fold Pfister form $\varphi$,
that is $q=q'\otimes\varphi$ for some form $q'$, 
there exists a direct summand
$N$ of the motive $M(Q_q)$ of the projective quadric $Q_q$ associated with $q$
such that
$$M(Q_q)\simeq
\begin{cases}
N\otimes M(\pl^{2^m-1}_F), &\text{if }\dim q'\text{ is even};\\
(N\otimes M(\pl^{2^m-1}_F))\oplus
M(Q_{\varphi})\{\tfrac{\dim q}{2}-2^{m-1}\},&\text{if }\dim q'\text{ is odd}.
\end{cases}
$$
In view of the Milnor-Husemoller theorem mentioned in the beginning,
formula~\eqref{mquad}
implies a shortened proof of Vishik's
result for $m=1$.

\bibliographystyle{chicago}

\paragraph{\rm N.~Semenov, Mathematisches Institut der LMU M\"unchen}

\paragraph{\rm K.~Zainoulline, Mathematisches Institut der LMU M\"unchen, 
Theresienstr.~39, 80333 M\"unchen; e-mail: kirill@mathematik.uni-muenchen.de}

\end{document}